\documentclass[11pt,envcountsame,envcountsect,leqno]{svmult}
\usepackage{amsmath,amssymb}

\newcommand{\Z}{{\mathbb{Z}}}
\newcommand{\N}{{\mathbb{N}}}

\newcommand{\R}{{\mathbb{R}}}
\newcommand{\bT}{{\mathbb{T}}}

\newcommand{\lamb}{\mathbf{\lambda}}
\newcommand{\cb}{{\mathbf{c}}}
\newcommand{\la}{\lambda}
\newcommand{\cD}{\mathcal{D}}
\newcommand{\cH}{{\mathcal{H}}}

\newcommand{\cP}{{\mathcal{P}}}
\newcommand{\cL}{{\mathcal{L}}}
\newcommand{\cR}{{\mathcal{R}}}
\newcommand{\cLR}{{\mathcal{LR}}}
\newcommand{\cW}{{\mathcal{W}}}
\newcommand{\fL}{{\mathfrak{L}}}

\newcommand{\fR}{{\mathfrak{R}}}
\newcommand{\fS}{{\mathfrak{S}}}

\newcommand{\nega}{\operatorname{neg}}
\newcommand{\sh}{\operatorname{sh}}
\newcommand{\tla}{\widetilde{\lambda}}
\newcommand{\dom}{\operatorname{dom}}

\renewcommand{\leq}{\leqslant}
\renewcommand{\geq}{\geqslant}

\newtheorem{conjA}{Conjecture}

\newtheorem{conjAp}{Conjecture}

\begin{document}
\title*{On domino insertion and Kazhdan--Lusztig cells 
in type $B_n$} 
\titlerunning{Domino insertion and Kazhdan--Lusztig cells}

\author{C\'edric Bonnaf\'e\inst{1}\and Meinolf Geck\inst{2}\and
Lacrimioara Iancu \inst{2}\and Thomas Lam\inst{3}}

\institute{Universit\'e de Franche-Comt\'e,  UFR Sciences et
Techniques, 16 route de Gray, 25 030 Besan\c{c}on, France
\texttt{bonnafe@math.univ-fcomte.fr} \and
Department of Mathematical Sciences, King's College,
Aberdeen University, Aberdeen AB24 3UE, Scotland, U.K.
\texttt{m.geck@maths.abdn.ac.uk,l.iancu@maths.abdn.ac.uk}
\and Department of Mathematics, Harvard University, Cambridge, MA 02138, USA
\texttt{tfylam@math.harvard.edu}}
\authorrunning{Bonnaf\'e, Geck, Iancu and Lam}

\date{\today}

\maketitle


\begin{abstract}
Based on empirical evidence obtained using the {\sf CHEVIE} computer algebra
system, we present a series of conjectures concerning the combinatorial
description of the Kazhdan--Lusztig cells for type $B_n$ with unequal
parameters. These conjectures form a far-reaching extension of the
results of Bonnaf\'e and Iancu obtained earlier in the so-called ``asymptotic
case''. We give some partial results in support of our conjectures.
\end{abstract}

\section{Introduction and the main conjectures} \label{sec1}

Let $W$ be a Coxeter group, $\Gamma$ be a totally ordered abelian group and
$L\colon W \rightarrow \Gamma_{\geq 0}$ be a
weight function, in the sense of Lusztig \cite[\S 3.1]{Lusztig03}.
This gives rise
to various pre-order relations on $W$, usually denoted by $\leq_{\cL}$,
$\leq_{\cR}$ and $\leq_{\cLR}$. Let $\sim_{\cL}$, $\sim_{\cR}$ and
$\sim_{\cLR}$ be the corresponding equivalence relations. The equivalence
classes are called the left, right and two-sided cells of $W$,
respectively. They were first defined by Kazhdan and Lusztig \cite{KaLu} in
the case where $L$ is the length function on $W$ (the ``equal parameter
case''), and by Lusztig \cite{Lusztig83} in general. They play a fundamental
role, for example, in the representation theory of finite or $p$-adic
groups of Lie type; see the survey in \cite[Chap.~0]{Lusztig03}.

Our aim is to understand the dependence of the Kazhdan--Lusztig cells on the
weight function $L$. We shall be interested in the case where $W$ is a finite
Coxeter group. Then unequal parameters can only arise in type $I_2(m)$
(dihedral), $F_4$ or $B_n$. Now types $I_2(m)$ and $F_4$ can be dealt with
by computational methods; see \cite{my04}. Thus, as far as finite Coxeter
groups are concerned, the real issue is to study type $B_n$ with unequal
parameters. And in any case, this is the most important case with respect
to applications to finite classical groups (unitary, symplectic, and
orthogonal). Quite recently, new connections between Kazhdan--Lusztig
cells in type $B_n$ and the theory of rational Cherednik algebras appeared
in the work of Gordon and Martino \cite{GM}.

The purpose of this paper is to present a series of conjectures
which would completely and explicitly determine the Kazhdan--Lusztig
cells in type $B_n$ for any positive weight function $L$. We will
also establish some relative results in support of these
conjectures. So let now $W=W_n$ be a Coxeter group of type $B_n$,
with generating set $S_n=\{t,s_1,\ldots,s_{n-1}\}$ and Dynkin
diagram as given below; the ``weights'' $a,b\in \Gamma_{> 0}$ attached
to the generators of $W_n$ uniquely determine a weight function
$L=L_{a,b}$ on $W_n$.
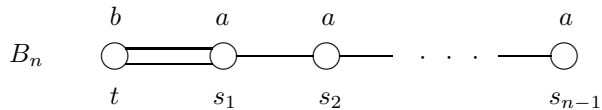
\begin{center}
\begin{picture}(250,42)
\put(  0, 18){$B_n$}
\put( 40, 20){\circle{10}}
\put( 44, 17){\line(1,0){33}}
\put( 44, 23){\line(1,0){33}}
\put( 81, 20){\circle{10}}
\put( 86, 20){\line(1,0){29}}
\put(120, 20){\circle{10}}
\put(125, 20){\line(1,0){20}}
\put(155, 17){$\cdot$}
\put(165, 17){$\cdot$}
\put(175, 17){$\cdot$}
\put(185, 20){\line(1,0){20}}
\put(210, 20){\circle{10}}
\put( 38, 32){$b$}
\put( 78, 32){$a$}
\put(118, 32){$a$}
\put(208, 32){$a$}
\put( 38, 02){$t$}
\put( 77, 02){$s_1$}
\put(117, 02){$s_2$}
\put(204, 02){$s_{n{-}1}$}
\end{picture}
\end{center}
If $b$ is ``large'' with respect to $a$, more precisely, if $b > (n-1) a$, 
then we are in the ``asymptotic case'' studied in \cite{BI} (see also 
\cite[Prop. 5.1 and Cor. 5.2]{BI2} for the determination of the exact bound). 
In general, we expect that the combinatorics governing the cells in type 
$B_n$  are provided by the
\begin{center}
``{\em domino insertion of a signed permutation into a $2$-core}'';
\end{center}
see \cite{Lam1}, \cite{vL}, \cite{ShWh} (see also \S\ref{sec3}).
Having fixed $r\geq 0$, let $\delta_r$ be the partition with parts
$(r,r-1, \ldots,0)$ (a $2$-core). Let $\cP_r(n)$ be the set of
partitions $\lambda \vdash (\frac{1}{2}r(r+1)+2n)$ such that
$\lambda$ has $2$-core $\delta_r$. Then the domino insertion with
respect to $\delta_r$ gives a bijection from $W_n$ onto the set of
all pairs of standard domino tableaux of the same shape $\lambda \in
\cP_r(n)$. We write this bijection as $w \rightarrow
(P^r(w),Q^r(w))$; see \cite[\S 2]{Lam1} for a detailed description.

The following conjectures have been verified for $n \leq 6$ by explicit 
computation using {\sf CHEVIE} \cite{chv} and the program {\sf Coxeter}
developed by du Cloux \cite{fokko}. For the basic definitions  
concerning Kazhdan--Lusztig cells, see Lusztig \cite{Lusztig03}.

\begin{conjA} \label{conj1}
Assume that $\Gamma=\Z$, $a=2$ and $b=2r+1$ where $r \geq 0$. Then the 
following hold.
\begin{itemize}
\item[(a)] $w,w'\in W_n$ lie in the same Kazhdan--Lusztig left cell if and
only if $Q^r(w)=Q^r(w')$.
\item[(b)] $w,w'\in W_n$ lie in the same Kazhdan--Lusztig right cell if and
only if $P^r(w)=P^r(w')$.
\item[(c)] $w,w'\in W_n$ lie in the same Kazhdan--Lusztig two-sided cell
if and only if all of $P^r(w)$, $Q^r(w)$, $P^r(w')$, $Q^r(w')$ have the
same shape.
\end{itemize}
\end{conjA}

\begin{remark} \label{rem1}
The $2$-core $\delta_r$, the set of partitions $\cP_r(n)$, and the
parameters $a=2$, $b=2r+1$ (where $\Gamma=\Z$) naturally arise in the
representation theory of the finite unitary groups $\mbox{GU}_N(q)$, 
where $N=\frac{1}{2}r(r+1)+2n$. The Hecke algebra of type $B_n$ with 
parameters $q^{2r+1},q^2,\ldots,q^2$ appears as the endomorphism algebra 
of a certain induced cuspidal representation. The irreducible 
representations of this endomorphism algebra parametrize the unipotent 
representations of $\mbox{GU}_N(q)$ indexed by partitions in $\cP_r(n)$; 
see \cite[\S 13.9]{Ca2}. In this case, Conjecture~A(c) 
is somewhat more precise than \cite[Conj. 25.3 (b)]{Lusztig03} 
(see \S\ref{sec:lusztig} for more details).
\end{remark}

\begin{conjAp} \label{conjA+}
Let $r \geq 0$ and assume that $a,b$ are any elements of $\Gamma_{>0}$
such that $ra<b<(r+1)a$. Then the statements in Conjecture~A still hold. 
That is, the Kazhdan-Lusztig (left, right, two-sided) cells for this choice 
of parameters coincide with those obtained for the special values $a=2$ 
and $b=2r+1$ (where $\Gamma=\Z$).
\end{conjAp}

\begin{remark}\label{order LR}
Assume we are in the setting of Conjecture~A or $\text{A}^+$. If $w \in W_n$, 
let $\lamb(w) \in \cP_r(n)$ denote the shape of $P^r(w)$ (or $Q^r(w)$). 
Let $\trianglelefteq$ denote the dominance order on partitions. The 
following property of the pre-order $\leq_{\cLR}$ has been checked for 
$n \leq 4$ by using {\sf CHEVIE} \cite{chv}:
\begin{equation*}
w \leq_{\cLR} w' \quad \mbox{if and only if} \quad  \lamb(w) 
\trianglelefteq \lamb(w') \tag{c$^+$}
\end{equation*}
\end{remark}

\begin{remark} \label{two?}
Assume that the statement concerning the left cells in Conjecture~A
(or $\text{A}^+$) is true. Since $P^r(w^{-1})=Q^r(w)$ (see for instance 
\cite[Lemma 7]{Lam1}), this would imply that the statement concerning 
the right cells is also true. However, it is not clear that the partition 
into two-sided cells easily follows from the knowledge of the partitions 
into left and right cells. Indeed, it is conjectured (but not proved in 
general) that the relation $\sim_{\cLR}$ is generated by $\sim_\cL$ and 
$\sim_\cR$.  This would follow from Lusztig's Conjectures (P4), (P9), 
(P10) and (P11).
\end{remark}

\begin{remark} \label{rem2}
If $b>(n-1)a$ (``asymptotic case''), then domino insertion is equivalent to
the generalized Robinson--Schensted correspondence in \cite[\S 3]{BI}
(see Theorem \ref{thm:BI}). Thus, Conjectures~A and $\text{A}^+$ holds in this 
case \cite[Th. 7.7]{BI}, \cite[Cor. 3.6 and Rem. 3.7]{BI2}. Also, the 
refinement (c$^+$) proposed in Remark~\ref{order LR} holds in this case 
if $w$ and $w'$ have the same $t$-length \cite[Th. 3.5 and Rem. 3.7]{BI2} 
(the {\it $t$-length} of an element $w \in W_n$ is the number of occurrences 
of $t$ in a reduced decomposition of $w$).
\end{remark}

\begin{remark} \label{cellular}
Assume that Conjectures~A and $\text{A}^+$ hold. Then we also conjecture that
the Kazhdan--Lusztig basis of the Iwahori--Hecke algebra $\cH_n$ associated
to $W_n$ and the weight function $L_{a,b}$ is a {\em cellular basis} in 
the sense of Graham--Lehrer \cite{GrLe}. See Subsection~\ref{sub-cell} for 
a more precise statement and applications to the representation theory of 
non-semisimple specialisations of $\cH_n$.
\end{remark}

We define the equivalence relation $\simeq_r$ on elements of $W_n$ as follows:
we write $w \simeq_r w'$ if and only if $Q^r(w)=Q^r(w')$. An equivalence class 
for the relation $\simeq_r$ is called a {\it left $r$-cell}. In other
words, left $r$-cells are the fibers of the map $Q^r$. Similarly, we define
{\it right $r$-cells} as the fibers of the map $P^r$ and
{\it two-sided $r$-cell} as the fibers of the map $\lamb : W_n \to \cP_r(n)$.

Conjectures~A and $\text{A}^+$ deal with the Kazhdan--Lisztig cells
for parameters such that $ra<b<(r+1)a$. The next conjecture is 
concerned with the Kazhdan--Lusztig cells whenever $b \in \N^* a$.

\begin{conjA} \label{conj3}
Assume that $b=ra$ for some $r\geq 1$. Then the Kazhdan--Lusztig left (resp. 
right, resp. two-sided) cells of $W_n$ are the smallest subsets of $W_n$ 
which are at the same time unions of left (resp. right, resp. two-sided) 
$(r-1)$-cells and left (resp. right, resp. two-sided) $r$-cells.
\end{conjA}

We will give a combinatorially more precise version of
Conjecture~B in \S\ref{sec:cycles}.

\begin{remark} \label{rem3}
(a) If $r \geq n$ then, since the left $r$-cells and 
the left $(r-1)$-cells coincide, then the Conjecture~B holds 
(``asymptotic case'', see Remark \ref{rem2}).

(b) There is one case which is not covered by Conjectures~A, $\text{A}^+$ 
or~B: it is when $b > ra$ for every $r \in \N$. But this case
is exactly the case which is dealt with in \cite[Th. 7.7]{BI} (and
\cite[Cor. 3.6]{BI2} for
the determination of two-sided cells) and it leads
to the same partition into left and two-sided cells as the case where
$(a,b) = (2, 2n-1)$ for instance (see Remark \ref{rem2}).

(c) The fundamental difference between the cases where $b \in \{a,2a,
\dots,(n-1)a \}$ and $b \not\in \{a,2a,\dots, (n-1) a\}$ is already 
appearant in \cite[Chap.~22]{Lusztig03}, where the ``constructible
representations'' are considered. Conjecturally, these are precisely
the representations given by the various left cells of $W$. By 
\cite[Chap. 22]{Lusztig03}, the constructible representations are all 
irreducible if and only if $b \not\in \{a,2a,\dots, (n-1) a\}$.

(d) Again, in Conjecture B, the statement concerning
left cells is equivalent to the statement concerning right cells.
However, the statement concerning two-sided cells would then follow
if one could prove that the relation $\sim_{\cLR}$ is generated
by the relations $\sim_\cL$ and $\sim_\cR$.

(e) Conjectures~$\text{A}^+$ and B are consistent with analogous 
results for type $F_4$ (see \cite{my04} as far as Conjecture~$\text{A}^+$ is 
concerned; Geck also checked that an analogue of Conjecture~B holds 
in type $F_4$).
\end{remark}

In Section~2, we will discuss representation-theoretic issues related
to Conjecture~A. In Sections~3 and~4, we will present a number 
of partial results in support of our conjectures. 

\section{Leading matrix coefficients and cellular bases} \label{sec2}

Let $W$ be a finite Coxeter group with generating set $S$.
Let $\Gamma$ be a totally ordered abelian group.
Let $L \colon W \rightarrow \Gamma$ be a weight function in the sense of
Lusztig \cite[\S 3.1]{Lusztig03}. Thus, we have $L(ww')=L(w)+L(w')$ for all
$w,w'\in W$ such that $l(ww')=l(w)+l(w')$ where $l\colon W \rightarrow \N$ is
the usual length function with respect to $S$ (where $\N=\{0,1,2, \ldots\}$).
Let $A={\Z}[\Gamma]$ be the group ring of $\Gamma$. It will be denoted 
exponentially: in other words, $A=\oplus_{\gamma \in \Gamma} \Z v^\gamma$ and 
$v^\gamma v^{\gamma'}=v^{\gamma+\gamma'}$.  If $\gamma_0 \in \Gamma$, let 
$A_{> \gamma_0} = \oplus_{\gamma > \gamma_0} \Z v^\gamma$.  We define 
similarly $A_{\geq \gamma_0}$, $A_{< \gamma_0}$ and $A_{\leq \gamma_0}$.

Let $\cH=\cH_A(W,S,L)$ be the corresponding
Iwahori--Hecke algebra. Then $\cH$ is free over $A$ with basis
$(T_w)_{w \in W}$; the multiplication is given by the rule
\[ T_sT_w=\left\{\begin{array}{cl} T_{sw} & \quad \mbox{if $l(sw)=l(w)+1$},
\\ T_{sw}+(v^{L(s)}-v^{-L(s)})T_w & \quad \mbox{if $l(sw)=l(w)-1$},\end{array}
\right.\]
where $w\in W$ and $s\in S$.  For basic properties of $W$ and $\cH$, we
refer to \cite{ourbuch}.

\subsection{Leading matrix coefficients}
We now recall the basic facts concerning the leading matrix coefficients
introduced in \cite{my02}. First, since $\Gamma$ is an ordered group, the 
ring $A$ is integral. Similarly, the group algebra $\R[\Gamma]$ is integral; 
we denote by $K=\R(\Gamma)$ its field of fractions.

Extending scalars from $A$ to the field $K$, we obtain a finite dimensional
$K$-algebra $\cH_K=K\otimes_A \cH$, with basis $(T_w)_{w \in W}$. It is
well-known that $\cH_K$ is split semisimple and abstractly isomorphic to the
group algebra of $W$ over $K$; see, for example, \cite[Remark 3.1]{GI}. Let
$\mbox{Irr}(\cH_K)$ be the set of irreducible characters of $\cH_K$. We write
this set in the form
\[ \mbox{Irr}(\cH_K)=\{\chi_\lambda \mid \lambda\in \Lambda\},\]
where $\Lambda$ is some finite indexing set. If $\lambda \in \Lambda$, we 
denote by $d_\lambda$ the degree of $\chi_\lambda$. We have a symmetrizing 
trace $\tau \colon \cH_K \rightarrow K$ defined by $\tau(T_1)=1$ and 
$\tau(T_w)=0$ for $1\neq w\in W$; see \cite[\S 8.1]{ourbuch}. The fact that 
$\cH_K$ is split semisimple yields that
\[ \tau=\sum_{\lambda\in \Lambda} \frac{1}{c_\lambda} \, \chi_\lambda
\qquad \mbox{where $0 \neq c_\lambda \in {\R}[\Gamma]$}.\]
The elements $c_\lambda$ are called the {\em Schur elements}. There is a
unique $a(\lambda)\in \Gamma_{\geq 0}$ and a positive real number $r_\lambda$ 
such that
\[ c_\lambda \in r_\lambda \, v^{-2a(\lambda)}+A_{> -2a(\lambda)};\]
see \cite[Def.~3.3]{my02}. The number $a(\lambda)$ is called the
$a$-invariant of $\chi_\lambda$. Using the {\em orthogonal representations}
defined in \cite[\S 4]{my02}, we obtain the {\em leading matrix coefficients}
$c_{w,\lambda}^{ij}\in \R$ for $\lambda\in \Lambda$ and $1\leq i,j
\leq d_\lambda$. See \cite[\S 4]{my02} for further general results
concerning these coefficients.

Following \cite[Def.~3.3]{GI}, we say that
\begin{itemize}
\item $\cH$ is {\em integral} if $c_{w,\lambda}^{ij}\in \Z$ for all
$\lambda\in \Lambda$ and $1\leq i,j\leq d_\lambda$;
\item $\cH$ is {\em normalized} if $r_\lambda=1$ for all $\lambda\in \Lambda$.
\end{itemize}
The relevance of these notions is given by the following result.

\begin{theorem}[See \protect{\cite[\S 4]{my02}} and
\protect{\cite[Lemma~3.8]{GI}}] \label{thm1}
Assume that $\cH$ is integral and normalized.
\begin{itemize}
\item[(a)] We have $c_{w,\lambda}^{ij}\in \{0,\pm 1\}$ for all $w\in W$,
$\lambda \in \Lambda$ and $1\leq i,j \leq d_\lambda$.
\item[(b)] For any $\lambda\in\Lambda$ and $1\leq i,j \leq d_\lambda$,
there exists a unique $w\in W$ such that $c_{w,\lambda}^{ij} \neq 0$;
we denote that element by $w=w_\lambda(i,j)$. The correspondence
$(\lambda,i,j)\mapsto w_\lambda(i,j)$ defines a bijective map
\[ \{(\lambda,i,j)\} \mid \lambda \in \Lambda,1\leq i,j \leq d_\lambda\}
\longrightarrow W.\]
\item[(c)] For a fixed $\lambda \in \Lambda$ and $1\leq k \leq d_\lambda$,
\begin{itemize}
\item[(i)] $\fL_{\lambda,k}:=\{w_\lambda(i,k) \mid 1 \leq i
\leq d_\lambda\}$ is contained in a left cell;
\item[(ii)] $\fR_{\lambda,k}:=\{w_\lambda(k,j) \mid 1 \leq j
\leq d_\lambda\}$ is contained in a right cell.
\end{itemize}
\end{itemize}
\end{theorem}


\begin{remark} \label{remRS} Assume that Lusztig's conjectures
{\bf (P1)--(P15)} in \cite[\S 14.2]{Lusztig03} hold for $\cH$.
Assume also that $\cH$ is normalized and integral. Combining
\cite[Corollary~4.8]{my05a} and \cite[Lemma~3.10]{GI}, we conclude that
the sets $\fL_{\lambda,k}$ and $\fR_{\lambda,k}$ are precisely the left
cells and the right cells of $W$, respectively.
\end{remark}


Now let $W=W_n$ be the Coxeter group of type $B_n$ as in Section~\ref{sec1};
let $\cH_n$ be the associated Iwahori--Hecke algebra with respect to the
weight function $L=L_{a,b}$ where $a,b\geq 0$.

\begin{proposition} \label{prop1}
Assume that $a>0$ and $b\not\in \{a,2a,\dots,(n-1)a\}$. Then $\cH_n$ is 
integral and normalized.
\end{proposition}

\begin{proof} The fact that $\cH_n$ is normalized follows from the explicit
description of $a(\lambda)$ in \cite[Prop.~22.14]{Lusztig03}. To show that 
$\cH_n$ is integral we follow once more the discussion in 
\cite[Example~3.6]{GI} where we showed that $\cH_n$ is integral if 
$b > (n-1) a$. So we may, and we will, assume from now on that $b < (n-1) a$.
Since $b \not\in \{a,2a,\dots,(n-1)a\}$, there exists a unique $r \geq 0$ 
such that $ra < b < (r+1)a$. Given $\lambda \in \Lambda$, let 
$\tilde{S}^\lambda$ be the Specht module constructed by Dipper--James--Murphy 
\cite{DJM3}. There is a non-degenerate $\cH_n$-invariant bilinear form 
$\langle \;,\; \rangle_\lambda$ on $\tilde{S}^\lambda$. Let $\{f_t \mid t 
\in \bT_\lambda\}$ be the orthogonal basis constructed in 
\cite[Theorem~8.11]{DJM3}, where $\bT_\lambda$ is the set of all standard 
bitableaux of shape $\lambda$.  Using the recursion
formula in \cite[Prop.~3.8]{djm2}, it is straightforward to show that, for
each basis element $f_t$, there exist integers $s_t, a_{ti}, b_{tj},c_{tk},
d_{tl} \in \Z$ such that $a_{ti} \geq 0$,  $b_{tj} \geq 0$, and
\[ \langle f_t,f_t\rangle_\lambda=v^{2s_ta}\cdot \frac{\prod_i
(1+v^{2a}+ \cdots +v^{2a_{ti}a})}{\prod_j (1+v^{2a}+\cdots +v^{2b_{tj}a})}
\cdot \frac{\prod_k \bigl(1+v^{2(b+c_{tk}a)}\bigr)}{\prod_l\bigl(1+
v^{2(b+d_{tl}a)}\bigr)}.\]
In \cite[Example~3.6]{GI}, we noticed that we also have $b+c_{tk}a>0$
and $b+d_{tl}a>0$ if $b>(n-1)a$, and this allowed us to deduce that
$\cH_n$ is integral in that case. Now, if we only assume that
$ra < b < (r+1) a$,
then $b+c_{tk}a$ and $b+d_{tl}a$ will no longer be strictly positive, but
at least we know that they cannot be zero. Thus, there exist $h_t,h_t',
m_{tk},m_{tl}'\in \Z$ such that
\begin{align*}
\prod_k \bigl(1+v^{2(b+c_{tk}a)}\bigr)&=v^{2h_t}\prod_k \bigl(1+v^{2m_{tk}}
\bigr) \qquad \mbox{where $m_{tk}>0$},\\
\prod_l \bigl(1+v^{2(b+d_{tl}a)}\bigr)&=v^{2h_t'}\prod_l \bigl(1+v^{2m_{tl}'}
\bigr) \qquad \mbox{where $m_{tl}'>0$}.
\end{align*}
Hence, setting
\[ \tilde{f}_t:=v^{-s_ta-h_t+h_t'} \cdot \Bigl(\prod_j (1+v^{2a}+ \cdots
+v^{2b_{tj}a})\Bigr)\cdot \Bigl(\prod_l(1+v^{2m_{tl}'})\Bigr)\cdot f_t,\]
we obtain $\langle \tilde{f}_t,\tilde{f}_t\rangle_\lambda \in 1+v{\Z}[v]$
for all $t$. We can then proceed exactly as in \cite[Example~3.6]{GI} to
conclude that $\cH_n$ is integral. \qed
\end{proof}

The above result, in combination with Theorem~\ref{thm1}, provides a 
first approximation to the left and right cells of $W_n$. By 
Remark~\ref{remRS}, the sets $\cL_{\lambda,k}$ and $\cR_{\lambda,k}$ should
be precisely the left and right cells, respectively. In this context,
Conjecture~A would give an explicit combinatorial description of 
the correspondence $(\lambda,i,j) \mapsto w_\lambda(i,j)$.

\subsection{Cellular bases} \label{sub-cell}
Let us assume that we are in the setting of Conjecture~A. As 
announced in Remark~\ref{cellular}, we believe that then the Kazhdan--Lusztig 
basis of $\cH_n$ will be cellular in the sense of Graham--Lehrer \cite{GrLe}.
To state this more precisely, we have to introduce some further notation.
Let $(C_w)_{w\in W}$ be the Kazhdan--Lusztig basis of $\cH_n$; the element
$C_w$ is uniquely determined by the conditions that 
\[ \overline{C}_w=C_w \qquad \mbox{and}\qquad C_w \equiv T_w \;\bmod
\cH_{n,>0},\]
where $\cH_{n,>0}=\sum_{w \in W_n} A_{>0} T_w$ and the bar denotes the 
ring involution defined in \cite[Lemma~4.2]{Lusztig03}. Furthermore, let 
$* : \cH_n \to \cH_n$ be the unique anti-automorphism such that 
$T_w^*=T_{w^{-1}}$ for all $w \in W_n$. We also have $C_w^*=C_{w^{-1}}$
for any $w\in W_n$.

Now assume that $a>0$ and $b \not\in \{a,2a,\ldots,(n-1)a\}$. If 
$b < (n-1) a$, let $r \geq 0$ be such that $ra < b < (r+1) a$. If $b > 
(n-1) a$, let $r$ be any natural number greater than or equal to $n-1$. 

We set $\Lambda_r:=\cP_r(n)$ and consider the partial order on $\Lambda_r$
given by the dominance order $\trianglelefteq$ on partitions. For $\lambda 
\in \Lambda_r$, let $M_r(\lambda)$ denote the set of standard domino 
tableaux of shape $\lambda$. If $(S,T) \in M_r(\lambda) \times 
M_r(\lambda)$, let $C_r(S,T):=C_w$ where $(S,T)= (P^r(w),Q^r(w))$.

\begin{conjA} \label{cell} With the above notation, $(\Lambda_r,M_r,C_r,*)$ 
is a cell datum in the sense of Graham--Lehrer \cite[Def.~1.1]{GrLe}.
\end{conjA}

The existence of a cellular structure has strong representation-theoretic 
applications. For the remainder of this section, assume that 
Conjecture~C is true.  Let $\theta \colon A \rightarrow k$  be a 
ring homomorphism into a field $k$. Extending scalars from $A$ to $k$, 
we obtain a $k$-algebra $\cH_{n,k}:=k \otimes_A \cH_n$ which will no longer
be semisimple in general. The theory of cellular algebras \cite{GrLe} 
provides, for every $\lambda \in \Lambda_r$, a {\em cell module} $S^\lambda$ 
of $\cH_{n,k}$, endowed with an $\cH_{n,k}$-equivariant bilinear form 
$\phi^\lambda$. We set 
\[ D^\lambda:=S^\lambda/\mbox{rad}\, \phi^\lambda \qquad \mbox{for every
$\lambda \in \Lambda_r$}.\]
Let $\Lambda_r^\circ:=\{D^\lambda \mid \lambda \in \Lambda_r \mbox{ such 
that } \phi^\lambda\neq 0\}$. Then we have
\[ \mbox{Irr}(\cH_{n,k})=\{ D^\lambda \mid \lambda \in \Lambda_r^\circ\}; 
\qquad \mbox{see Graham--Lehrer \cite[Thm~3.4]{GrLe}}.\]
Thus, we obtain a natural parametrization of the irreducible representations
of $\cH_{n,k}$ by the set $\Lambda_r^\circ\subseteq \Lambda_r$.

\begin{remark} \label{cell1} Assume that $b>(n-1)a>0$. Then 
Conjecture~C holds by \cite[Cor.~6.4]{my05b}. In this 
case, the set $\Lambda_r^\circ$ is determined explicitly by 
Dipper--James--Murphy \cite{DJM3} and Ariki \cite{Ar2}. Finally, 
Iancu--Pallikaros \cite{IaPa} show that  the cell modules 
$S^\lambda$ are canonically isomorphic to the Specht modules 
defined by Dipper--James--Murphy \cite{DJM3}.
\end{remark}

\begin{remark} \label{cell2} 
Now consider arbitrary values of $a,b$ such that $a>0$ and 
$b \not\in \{a,2a, \ldots,(n-1)a\}$. Then, assuming that the 
conjectured relation (c$^+$) in Remark~\ref{order LR} holds, a 
description of the set $\Lambda_r$ follows from the results of 
Geck--Jacon \cite{GJ} on {\em canonical basic sets}. Indeed, one 
readily shows that the set $\Lambda_r^\circ$ coincides with the 
{\em canonical basic set} determined by \cite{GJ}. Thus, by the 
results of \cite{GJ}, we have explicit combinatorial 
descriptions of $\Lambda_r^\circ$ in all cases. Note that these 
descriptions heavily depend on $a,b$ and $\theta\colon 
A \rightarrow k$. 

It is shown in \cite{myprinc} that, if $a=2$ and $b=1$ or $3$, then the sets 
$\Lambda_r^\circ$ parametrize the modular principal series representations 
of the finite unitary groups.
\end{remark}

\section{Domino insertion} \label{sec3}

The aim of this section is to describe the domino insertion
algorithm and to provide some theoretical evidences for
Conjecture A. For this purpose we will see $W_n$
as the group of permutations $w$ of $\{-1,-2,\dots,-n\} \cup \{1,2,\dots,n\}$
such that $w(-i)=-w(i)$ for any $i$. The identification is as
follows: $t$ corresponds to the transposition $(1,-1)$ and
$s_i$ to $(i,i+1)(-i,-i-1)$. If $r \leq n$, we identify $W_r$
with the subgroup of $W_n$ generated by $S_r=\{t,s_1,s_2,\dots,s_{r-1}\}$.
The symmetric group of degree $n$ will be denoted by $\fS_n$: when necessary,
we shall identify it in the natural way with the subgroup of $W_n$
generated by $\{s_1,s_2,\dots,s_{n-1}\}$.
Let $t_1=t$ and, if $1 \leq i \leq n-1$, let $t_{i+1}=s_i t_i s_i$.
As a signed permutation, $t_i$ is just the transposition $(i,-i)$.

\begin{remark}\label{descent}
Since we shall be interested in various descent sets of elements of
$W_n$, we state here for our future needs the following two easy
facts. Let $w \in W_n$. Then the following hold.
\begin{itemize}
\item[(a)] If $1 \leq i \leq n-1$, then $\ell(ws_i) > \ell(w)$ if and only if
$w(i) < w(i+1)$.
\item[(b)] If $1 \leq i \leq n$, then $\ell(w t_i) > \ell(w)$ if and only
if $w(i) > 0$.
\end{itemize}
\end{remark}

\subsection{Partitions and Tableaux}
\label{sec:partitions} We refer to~\cite{Lam1,ShWh} for further
details of the material in this section.  We shall assume some
familiarity with (standard) Young tableaux.

Let $\la = (\la_1 \geq \la_2 \geq \ldots \geq \la_{l(\la)} > 0)$ be
a partition of $n = |\la| = \la_1 + \la_2 + \cdots + \la_{l(\la)}$. We will
not distinguish between a partition $\la$ and its {\it Young
diagram} (often denoted $D(\la)$).  Our Young diagrams will be drawn
in the English notation so that the boxes are upper-left justified.
When $\la$ and $\mu$ are partitions satisfying $\mu \subset \la$ we
will use $\la/\mu$ to denote the shape corresponding to the
set-difference of the diagrams of $\la$ and $\mu$. We call $\la/\mu$
a {\it domino} if it consists of exactly two squares sharing an
edge.

The {\it $2$-core} (or just core) $\tla$ of a shape $\la$ is
obtained by removing dominoes from $\la$, keeping the shape a
partition, until this is no longer possible.  The partition $\tla$
does not depend on how these dominoes are removed.  Every 2-core has
the shape of a staircase $\delta_r = (r,r-1,\ldots,0)$ for some
integer $r \geq 0$.

We denote the set of partitions by $\cP$ and the set of partitions
with 2-core $\delta_r$ by $\cP_r$.  The set of all partitions $\la$
satisfying the conditions:
$$
\tla = \delta_r \ \mbox{and} \ |\la| = |\delta_r| + 2n
$$
will be denoted $\cP_r(n)$.  Note that $\cP = \cup_{r,n} \cP_r(n)$ is a
disjoint union.

\medskip

A {\it (standard) domino tableau}  $D$ of {\it shape} $\la \in
\cP_r(n)$ consists of a tiling of the shape $\la / \tla$ by dominoes
and a filling of the dominoes with the integers $\{1,2,\ldots,n\}$,
each used exactly once, so that the numbers are increasing when read
along either the rows or columns.  The {\it value} of a domino is
the number written inside it.  We will denote by $\dom_i$ the domino
with the value $i$ inside. We will also write $\sh(D) = \la$ for the
shape of $D$. An equivalent description of the domino tableau $D$ is
as the sequence of partitions $\{\tla=\la^0 \subset \la^1 \subset
\ldots \subset \la^n = \la\}$, where $\sh(\dom_i) =
\la^i/\la^{i-1}$.  If the values of the dominoes in a tableau $D$
are not restricted to the set $\{1,2,\ldots,n\}$ (but each value
occurs at most once), we will call $D$ an {\it injective domino
tableau}.

\medskip

We now describe a number of operations on standard Young and domino
tableaux needed in the sequel.  One may obtain a standard Young
tableau $T = T(D)$ from a domino tableau $D$ by replacing a domino
with the value $i$ in $D$ by two boxes containing $\bar i$ and $i$
in $T$. The boxes are placed so that $T$ is standard with respect to
the order $\bar 1 < 1 < \bar 2 < 2 < \cdots$.  If $D$ has shape
$\lambda$ then $T(D)$ will have shape $\lambda/\tla$.  Suppose now
that $Y$ is a standard Young tableau of shape $\tla$ filled with
letters smaller than any of the letters occurring in $D$.  Define
$T_Y(D)$ by ``filling'' in the empty squares in $T(D)$ with the
tableau $Y$.

Let $T$ be a standard Young tableau and $i$ a letter occurring in
$T$.  The {\it conversion} process proceeds as follows
(see~\cite{Hai, ShWh}). Replace the letter $i$ in $T$ with another
letter $j$.  The resulting tableau may not be standard, so we
repeatedly swap $j$ with its neighbours until the tableau is
standard. We say that the value $i$ has been converted to $j$.

Now let $T$ be any standard Young tableau filled with barred $\bar
i$ and non-barred letters $i$.  Define $T^{\nega}$ by successively
converting barred letters $\bar i$ to negative letters $-i$,
starting with the smallest letters.  The main fact that we shall
need is that the operation ``$\nega$'' is invertible.  We refer the
reader to~\cite{ShWh} for a full discussion of these operations.



\subsection{The Barbasch-Vogan domino insertion algorithm}
The Robinson-Schensted correspondence establishes a bijection 
\[ \pi \leftrightarrow (P(\pi),Q(\pi))\]
between permutations $\pi \in
\fS_n$ and pairs of standard tableaux with the same shape and size
$n$ (see~\cite{Sta}). Domino insertion generalizes this by replacing
the symmetric group with the hyperoctahedral group.  It depends on
the choice of a core $\delta_r$, and establishes a bijection between
$W_n$ and pairs $(P^r,Q^r)$ of standard domino tableaux of the same
shape $\la \in \cP_r(n)$.  There are in fact many such bijections but
we will be concerned only with the algorithm introduced by Barbasch
and Vogan~\cite{BV} and later given a different description by
Garfinkle~\cite{gar1}.  We now describe this algorithm following the
more modern expositions~\cite{Lam1, ShWh}.

Let $D$ be an injective domino tableau with shape $\la$ such that $i
> 0$ is a value which does not occur in $D$.  We describe the insertion
$E = D \leftarrow i$ (or $E = D \leftarrow -i$) of a horizontal
(vertical) domino with value $i$ into $D$.  Let $D_{<i} \subset D$
denote the sub-domino tableau of $D$ containing all dominoes with
values less than $i$.  If $\la$ has a 2-core $\tilde{\la}$,
then we will always assume that $\tilde{\la} \subset \sh(D_{<i})$.
Let $E_{\leq i}$ be the domino tableau obtained from $D_{<i}$ by
adding an additional vertical domino in the first column or an
additional horizontal domino in the first row labeled $i$.

For $j > i$ we define $E_{\leq j}$, supposing that $E_{\leq j-1}$ is
known.  If $D$ contains no domino labeled $j$ then $E_{\leq
j} = E_{\leq j-1}$; otherwise let $\dom_j$ denote the domino in $D$
labeled $j$.  Let $\mu = \sh(E_{\leq j-1})$.  We now distinguish
four cases:
\begin{enumerate}
\item If $\mu \cap \dom_j = \emptyset$ do not touch, then we set $E_{\leq j} =
E_{\leq j-1} \cup \dom_j$.
\item If $\mu \cap \dom_j = (k,l)$ is exactly one
square in the $k$-th row and $l$-th column, then we add a domino
containing $j$ to $E_{\leq j-1}$ to obtain the tableau $E_{\leq j}$
which has shape $\mu \cup \dom_j \cup (k+1,l+1)$.
\item If $\mu \cap \dom_j = \dom_j$ and $\dom_j$ is horizontal, then we
bump the domino $\dom_j$ to the next row, by setting $E_{\leq j}$ to
be the union of $E_{\leq j-1}$ with an additional (horizontal)
domino with value $j$ one row below that of $\dom_j$.
\item If $\mu \cap \dom_j = \dom_j$ and $\dom_j$ is vertical, then we
bump the domino $\dom_j$ to the next column, by setting $E_{\leq j}$
to be the union of $E_{\leq j-1}$ with an additional (vertical)
domino with value $j$ one column to the right of $\dom_j$.
\end{enumerate}
Finally we let $E = \lim_{j \to \infty} E_{\leq j}$.

\medskip

Let $w = w(1) w(2) \cdots w(n) \in W_n$ be a hyperoctahedral
permutation written in one-line notation. Thus, for each $i$, we have
$w(i) \in \{\pm 1,\pm 2,\ldots,\pm n\}$; furthermore, $|w(1)||w(2)|
\cdots|w(n)| \in \fS_n$ is a usual permutation. Let $\delta_r$ be a 
$2$-core assumed to be
fixed. Then the insertion tableau $P^r(w)$ is defined as
$((\ldots((\delta_r \leftarrow w(1)) \leftarrow w(2)) \cdots )
\leftarrow w(n))$. The sequence of shapes obtained in the process
defines another standard domino tableau called the recording tableau
$Q^r(w)$ of $w \in W_n$. The insertion tableau $P^r(w)$ can of
course be defined for any sequence $w = w(1) w(2) \cdots w(n)$ such
that $|w(i)| \neq |w(j)|$ for $i \neq j$.

The following theorem is due to Barbasch-Vogan~\cite{BV} and
Garfinkle~\cite{gar1} when $r = 0,1$ and extended by van
Leeuwen~\cite{vL} to larger cores.

\begin{theorem}
\label{thm:bij} Fix $r \geq 0$.  The domino insertion algorithm
defines a bijection between $w \in W_n$ and pairs $(P,Q)$ of
standard domino tableaux of the same shape lying in $\cP_r(n)$. This
bijection satisfies the equality $P^r(w) = Q^r(w^{-1})$.
\end{theorem}

%

It is easy to see that the bijectivity in Theorem~\ref{thm:bij}
together with Conjecture~$\text{A}^+$ would imply that the relevant
left cell representations are irreducible. This is consistent with
the conjecture that left cell representations for $W_n$ are
irreducible for ``generic parameters'' and in particular if $b
\notin \{a,2a,\dots,(n-1)a\}$ (see Proposition~\ref{prop1}).

We have computational evidences for Conjectures A, $\text{A}^+$ and~B: 
they were checked for $n \leq 6$ by using {\sf CHEVIE} \cite{chv} and 
{\sf Coxeter} \cite{fokko}. In the rest of this section, we shall give 
theoretical evidences for Conjecture~A and $\text{A}^+$ (induction of 
cells, multiplication by the longest element, link to 
\cite[Conj. 25.3]{Lusztig03}, asymptotic case, quasi-split case, 
right descent sets, coplactic relations).

\subsection{Conjecture A and Lusztig's Conjecture 25.3}
\label{sec:lusztig} There is an alternative description (in the case
where $r=0$, $1$, it is in fact the original description of Barbasch
and Vogan) of domino insertion. As we will now explain, it is related
to \cite[Conj. 25.3]{Lusztig03}. Let us fix in this subsection a
Coxeter group $(W,S)$ of type $A_{2n+r(r+1)/2-1}$. Let $\sigma$ be the
unique non-trivial automorphism of $W$ such that $\sigma(S)=S$. If $J$
is a subset of $S$, we denote by $W_J$ the parabolic subgroup of $W$
generated by $J$ and let $w_J$ denote the longest element of $W_J$.

Let $I$ be the unique connected (when we view it as a subdiagram of
the Dynkin diagram of $(W,S)$) subset of $S$ of cardinality
$r(r+1)/2 -1$ (or $0$ if $r=0$) such that $\sigma(I)=I$:
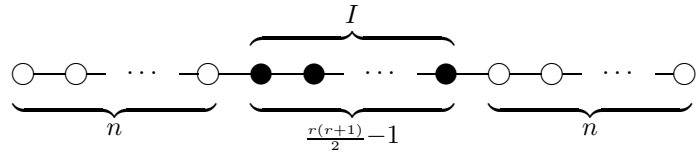
$$\begin{picture}(250,55)
\put(  0, 30){\circle{8}}
\put(  4, 30){\line(1,0){12}} 
\put( 20, 30){\circle{8}}
\put( 24, 30){\line(1,0){7}} 
\put( 39, 30){$\dots$}
\put( 70, 30){\circle{8}}
\put( 66, 30){\line(-1,0){7}}
\put( 74, 30){\line(1,0){12}} 
\put( 90, 30){\circle*{8}}
\put( 94, 30){\line(1,0){12}} 
\put(110, 30){\circle*{8}}
\put(114, 30){\line(1,0){7}} 
\put(129, 30){$\dots$} 
\put(160, 30){\circle*{8}}
\put(156, 30){\line(-1,0){7}}
\put(164, 30){\line(1,0){12}} 
\put(180, 30){\circle{8}}
\put(184, 30){\line(1,0){12}} 
\put(200, 30){\circle{8}}
\put(204, 30){\line(1,0){7}} 
\put(219, 30){$\dots$} 
\put(250, 30){\circle{8}}
\put(246, 30){\line(-1,0){7}} 
\put( -4, 20){$\underbrace{\hspace{2.7cm}}_{\displaystyle{n}}$} 
\put( 86, 20){$\underbrace{\hspace{2.7cm}}_{\frac{r(r+1)}{2}
\displaystyle{-1}}$}
\put(176, 20){$\underbrace{\hspace{2.7cm}}_{\displaystyle{n}}$}
\put( 86, 40){$\overbrace{\hspace{2.7cm}}^{\displaystyle{I}}$}
\end{picture}$$
Let $\cW$ denote the subgroup of $W$ consisting of all elements $w$ such
that $wW_Iw^{-1}=W_I$ and $w$ has minimal length in $w W_I$ (see
\cite[\S 25.1]{Lusztig03}). If $\Omega$ is a $\sigma$-orbit in $S\setminus
I$, we set $s_\Omega = w_{I \cup \Omega} w_I$. If $0 \leq i \leq n-1$, let
$\Omega_i$ denote the orbit of $\sigma$ in $S\setminus I$ consisting of
elements which are separated from $I$ by $i$ nodes in the Dynkin
diagram. Then $\{\Omega_0,\Omega_1,\dots,\Omega_{n-1}\}$ is the set of 
orbits of $\sigma$ in $S\setminus I$. Moreover, there is a unique morphism of
groups $\iota_r : W_n \to \cW^\sigma$ that sends $t$ to $s_{\Omega_0}$ and
$s_i$ to $s_{\Omega_i}$ (for $1 \leq i \leq n-1$). It is an isomorphism
of groups (see \cite[\S 25.1]{Lusztig03}).

The morphism $\iota_r$ can be described explicitly in the language
of signed permutations.  First identify $W$ with the permutation
group of the following $2n + r(r-1)/2$ elements (ordered according
to the ordering of $S$):
$$
\{-n < -(n-1) < \cdots < -1 < 0_1 < 0_2 < \cdots < 0_{r(r-1)/2} < 1
< 2 < \cdots n\}
$$
so that the subgroup $W_I$ (which is isomorphic to $\fS_{r(r+1)/2}$)
acts on the elements $\{0_1,0_2,\ldots,0_{r(r-1)/2}\}$.  Let $w =
w(1)w(2)\cdots w(n) \in W_n$.  Then the two-line notation of
$\iota_r(w)$ is given by \begin{equation}\label{eq:TLN}\left(
\begin{array}{ccccccccccc}
-n  & \cdots & -1 & 0_1 & \cdots & 0_{r(r-1)/2}& 1 & 2& \cdots &n
\\ -w(n) & \cdots & -w(1) & 0_1 & \cdots & 0_{r(r-1)/2}&
w(1) & w(2) & \cdots &w(n)
\end{array}\right).\end{equation}

Now, let $\cb_0$ denote the two-sided cell of $W_I$ which has
``shape'' $\delta_r$.  If $w$, $w' \in \fS_n$, we write $w \simeq_\fS
w'$ if $Q(w)=Q(w')$ (the equivalence relation $\simeq_\fS$ defines
the Robinson-Schensted left cells of $\fS_n$, which coincide with
the Kazhdan-Lusztig left cells \cite[\S 5]{KaLu}).

\begin{theorem}
\label{thm:sn} Fix $x \in \cb_0$. Let $w \in W_n$, $r \geq 0$ and
$\pi = \iota_r(w) x \in \fS_{2n+r(r+1)/2}$.  Then we have
$$
(T_{P(x)}(P^r(w)))^{\nega} = P(\pi) \hspace{15pt} \mbox{and}
\hspace{15pt} (T_{Q(x)}(Q^r(w)))^{\nega} = Q(\pi).
$$
Since $\nega$ is invertible, in particular $w \simeq_r w'$ if and
only if $\iota_r(w)x \simeq_{\fS} \iota_r(w') x$.
\end{theorem}

For the construction of $T_{P(x)}(P^r(w))$ in Theorem~\ref{thm:sn}
we are using the ordering $0_1 < 0_2 < \cdots < 0_{r(r-1)/2} < \bar
1 < 1 < \bar 2 < 2 < \cdots < \bar n < n$.  In the case where $r=0$
or $1$, Theorem~\ref{thm:sn} is essentially \cite[Theorem 4.2.3]{vL}
with different notation.  To generalize the result to all $r \geq 0$
we follow the approach of~\cite{ShWh}.

\begin{proof}
For the case $r = 0$ the theorem is exactly~\cite[Theorem~32]{ShWh}.
We now explain, assuming familiarity with~\cite{ShWh}, how to extend
the result to larger cores.  It is shown in~\cite[Lemma~31]{ShWh}
that a domino insertion $D \leftarrow i$ can be imitated by {\it
doubly mixed insertion}, denoted $P_{m^*}$. The proof
of~\cite[Lemma~31]{ShWh} is local, and remains valid when we replace
$D_{< i}$ by any Young tableau of the same shape, filled with
``small'' letters.  More precisely, their proof shows that
$T_{P(x)}(P^r(w))$ can be obtained by doubly mixed insertion of a
``biword'' $w^{\rm dup}$ (explicitly defined in~\cite{ShWh}) into
$P(x)$. Thus one has
\begin{equation}
\label{eq:DM} T_{P(x)}(P^r(w)) = P_{m^*}(x \sqcup w^{\rm dup}),
\end{equation}
where $a \sqcup b$ denotes the word obtained from concatenating $a$
and $b$.  In the notation of~\cite{ShWh}, $x$ here is a biword with
no bars so that $P_{m^*}(x) = P(x)$.

Now~\cite[Theorem 21 and Proposition 14]{ShWh} connect doubly mixed
insertion with usual Schensted insertion via the equation
\begin{equation}\label{eq:inin} 
P_{m^*}(u)^{\nega} = P(u^{{\rm inv} \nega {\rm inv} \nega}).
\end{equation} 
The operation denoted ``${\rm inv}$ $\nega$
${\rm inv}$ $\nega$'' in~\cite{ShWh} applied to $x \sqcup w^{\rm
dup}$ coincides with our inclusion $\iota_r(w)x$.  Combining
(\ref{eq:DM}) and (\ref{eq:inin}) one obtains
\[ (T_{P(x)}(P^r(w)))^{\nega} = P(\pi).\]  
The statement about recording tableau is obtained analogously, or by 
using the equation $Q(\pi) = P(\pi^{-1}) = P(x^{-1}\iota_r(w)^{-1}) =
P(\iota_r(w^{-1})x^{-1})$. \qed
\end{proof}

\begin{remark}\label{rem:independence}
Note that the last statement of Theorem~\ref{thm:sn} does not depend
on the choice of $x \in \cb_0$.
\end{remark}

%
%

\begin{corollary}\label{coro:25.3}
If $r \geq 0$ and if $(a,b)=(2,2r+1)$, then Conjecture A(c) 
agrees with \cite[Conj. 25.3]{Lusztig03} for the case
$(W,S,I,\sigma)$ described above.
\end{corollary}


\subsection{Longest element}
Let $w_0$ denote the longest element of $W_n$: it is equal to $t_1
t_2 \dots t_n$ (or to $-1$ $-2$ \dots $-n$ in the one line
notation). It is a classical result that two elements $x$ and $y$ in
$W_n$ satisfy $x \sim_\cL y$ if and only if $w_0 x \sim_\cL w_0 y$.
The next result shows that the relations $\simeq_r$ share the same
property.

\begin{proposition}\label{prop:longest}
Let $r \geq 0$ and let $x$, $y \in W_n$. Then $x \simeq_r y$ if and only
if $w_0 x \simeq_r w_0 y$.
\end{proposition}

\begin{proof}
This follows from the easy fact that $P^r(w_0 x)$ (resp. $Q^r(w_0
x)$) is the conjugate (that is, the transpose with respect to the
diagonal) of $P^r(x)$ (resp. $Q^r(x)$), and similarly for $y$. \qed
\end{proof}


\subsection{Induction of cells}
Let $m \leq n$. Let $X_m^n$ denote the set of elements $w \in W_n$ which
have minimal length in $w W_m$. It is a cross-section of $W_n/W_m$.
By Remark \ref{descent}, an element $x \in W_n$ belongs to $X_m^n$ if
and only if $0 < x(1) < x(2) < \dots < x(m)$.
A theorem of Geck \cite{geck induction} asserts that, if $C$ is a 
Kazhdan-Lusztig left cell of $W_m$ (associated with the restriction of 
$L_{a,b}$ to $W_m$), then $X_m^n C$ is a union of Kazhdan-Lusztig left
cells of $W_n$. The next result show that the same hold if we replace
{\it Kazhdan Lusztig left cell} by {\it left $r$-cell}.

\begin{proposition}\label{prop:induction}
Let $r \geq 0$. If $C$ is a left $r$-cell of $W_m$,
then $X_m^n C$ is a union of left $r$-cells of $W_n$.
\end{proposition}

\begin{proof}
Let $w$, $w' \in W_m$ and $x$, $x' \in X_m^n$ be such that
$xw \simeq_r x'w'$ (in $W_n$). We must show that
$w \simeq_r w'$ (in $W_m$). For the purpose of this proof,
we shall denote by $(P_n^r(w),Q_n^r(w))$ (resp. $(P_m^r(w),Q_m^r(w))$)
the pair of standard domino tableaux obtained by viewing $w$ as an element 
of $W_n$ (resp. of $W_m$).
Then, since $x$ is increasing on $\{1,2,\dots,m\}$ and takes only
positive values, the dominoes filled with $\{1,2,\dots,m\}$
in the recording tableau $Q_n^r(xw)$ are the same as the one
in the recording tableau $Q_n^r(w)$. In particular,
$Q_m^r(w)$ is obtained from $Q_n^r(xw)$ by removing the dominoes
filled by $\{m+1,m+2,\dots,n\}$. Similarly,
$Q_m^r(w')$ is obtained from $Q_n^r(x'w')$ by removing the dominoes
filled by $\{m+1,m+2,\dots,n\}$. Since $Q_n^r(xw)=Q_n^r(x'w')$
by hypothesis, we have that $Q_m^r(w)=Q_m^r(w')$. In other
words, $w \simeq_r w'$ in $W_m$. \qed
\end{proof}

\begin{corollary}\label{coro:parabolic}
Let $r \geq 0$ and let $x$ and $y$ be two elements of $W_m$.
Then $x \simeq_r y$ in $W_m$ if and only if $x \simeq_r y$ in $W_n$.
\end{corollary}

The previous corollary shows that it is not necessary to make the ambient
group precise when one studies the equivalence relation $\simeq_r$.

Geck's result \cite{geck induction} is valid for any Coxeter group
and any parabolic subgroup. We shall investigate now the analogue of
Proposition \ref{prop:induction} for the parabolic subgroup $\fS_n$
of $W_n$.  We denote by $X(n)$ the set of elements $w \in W_n$ which
have minimal length in $w\fS_n$. It is a cross-section of
$W_n/\fS_n$. By Remark \ref{descent}, an element $w \in W_n$ belongs
to $X(n)$ if and only if $w(1) < w(2) < \dots < w(n)$.

\begin{proposition}\label{prop:induction sn}
Let $r \geq 0$ and let $C$ be a Robinson-Schensted left cell of $\fS_n$.
Then $X(n) C$ is a union of domino left cells for $\simeq_r$.
\end{proposition}

\begin{proof}
Let $w, w' \in \fS_n$ and $x, x' \in X(n)$ be such that $xw \simeq_r
x'w'$ in $W_n$.  We must show that $w \simeq_\fS w'$.  It is well
known that for two words $a_1a_2\cdots a_k$ and $b_1b_2 \cdots b_k$
one has $Q(a_1a_2\cdots a_k) = Q(b_1 b_2 \cdots b_k) \implies Q(a_j
a_{j+1} \cdots a_l) = Q(b_j b_{j+1} \cdots b_l)$ for any $1 \leq j
\leq l \leq k$. Indeed this is Geck's result~\cite{geck induction}
for $\fS_n$.

By Theorem~\ref{thm:sn}, we have $Q(\iota_r(xw)c) = 
Q(\iota_r(x'w')c)$ for any $c \in \cb_0$.  
Treating $i_r(xw)c$ as a word using (\ref{eq:TLN}), we thus have
$$Q(xw(1)\;xw(2)\;\cdots\;xw(n)) =
Q(x'w'(1)\;x'w'(2)\;\cdots\;x'w'(n)).$$  But $x(1) < x(2) < \dots <
x(n)$ so that this is equivalent to $Q(w) = Q(w')$. \qed
\end{proof}

If $x \in X_m^n$, and if $w$, $w' \in W_m$ are such that 
$w \sim_\cL w'$ in $W_m$
then a result of Lusztig \cite[Prop. 9.13]{Lusztig03} asserts that
$wx^{-1} \sim_\cL w'x^{-1}$. The next result shows that the same
statement holds if we replace $\sim_\cL$ by $\simeq_r$.

\begin{proposition}\label{prop:domino lusztig}
Let $r \geq 0$, $x \in X_m^n$ and $w$, $w' \in W_m$ be such that
$w \simeq_r w'$. Then $wx^{-1} \simeq_r w'x^{-1}$.
\end{proposition}

\begin{proof}
Let us use here the notation of the proof of Proposition
\ref{prop:induction}. So we assume that
$P_m^r(w^{-1})=P_m^r(w^{\prime-1})$ and we must show that
$P_n^r(xw^{-1})=P_n^r(xw^{\prime -1})$. Let $D=((\ldots((\delta_r
\leftarrow xw^{-1}(1))\leftarrow xw^{-1}(2)) \cdots ) \leftarrow
xw^{-1}(m))$ and $D'=((\ldots((\delta_r \leftarrow
xw^{\prime-1}(1))\leftarrow xw^{\prime-1}(2)) \cdots ) \leftarrow
xw^{\prime-1}(m))$. Since $w^{-1}(i)=i$ if $i \geq m+1$, we have
$P_n^r(xw^{-1})=((\ldots((D \leftarrow x(m+1))\leftarrow x(m+2))
\cdots ) \leftarrow x(n))$ and similarly for $P_n^r(xw^{\prime
-1})$. Therefore, we only need to show that $D=D'$. But, since
$w^{-1}$ stabilizes $\{1,2,\dots,m\}$ and since $x$ is increasing on
$\{1,2,\dots,m\}$ and takes only positive values, it follows from
the domino insertion algorithm that $D$ is obtained from
$P_m^r(w^{-1})$ by applying $x$, or in other words replacing the
domino $\dom_i$ by $\dom_{x(i)}$. Similarly, $D'$ is obtained from
$P_m^r(w^{\prime -1})$ by applying $x$. Since
$P_m^r(w^{-1})=P_m^r(w^{\prime-1})$ by hypothesis, we get that
$D=D'$, as desired. \qed
\end{proof}

As for Geck's result, Lusztig's result \cite[Prop. 9.13]{Lusztig03} is valid
for any parabolic subgroup of any Coxeter group. The next result
is the analogue of Proposition \ref{prop:domino lusztig}
for the parabolic subgroup $\fS_n$ of $W_n$.

\begin{proposition}\label{prop:domino lusztig sn}
Let $r \geq 0$, $x \in X(n)$ and $w$, $w' \in \fS_n$ be such that
$w \simeq_\fS w'$. Then $wx^{-1} \simeq_r w'x^{-1}$.
\end{proposition}

\begin{proof}
We must show that $P^r(xw^{-1}) = P^r(xw^{\prime-1})$ knowing that
$P(w^{-1}) = P(w^{\prime-1})$.  For $u \in W_n$, denote by $u_{\fS}$
the word $u(1)\;u(2)\;\cdots\;u(n)$ and $u_{-\fS}$ the word $-u(n)
\; -u(n-1) \; \cdots\;-u(1)$.  The equation $P(w^{-1}) =
P(w^{\prime-1})$ gives
\begin{equation} \label{eq:PSG}
P((xw')_\fS) = P((xw^{\prime-1})_\fS) \hspace{5pt} \mbox{and}
\hspace{5pt} P((xw')_{-\fS}) = P((xw^{\prime-1})_{-\fS}),
\end{equation} where we are comparing pairs of standard Young
tableaux filled with the set of letters $\{x(1),x(2),\ldots,x(n)\}$
(resp. $(\{-x(1),-x(2),\ldots,-x(n)\}$).

By Theorem~\ref{thm:sn}, it is enough to show that
$P(\iota_r(xw^{-1})c)= P(\iota_r(xw^{\prime-1})c)$ for some fixed $c
\in \cb_0$. Using~(\ref{eq:TLN}), we may write $\iota_r(xw^{-1})c$ in
one-line notation as the concatenation $(xw^{-1})_{-\fS} \sqcup c
\sqcup (xw^{-1})_{\fS}$.  It is well known that if $a, a', b$ are
words such that $P(a) = P(a')$ then one has $P(a \sqcup b) = P(a'
\sqcup b)$ (this is Lusztig's result \cite[Prop. 9.13]{Lusztig03}
for the symmetric group).  Combining this with (\ref{eq:PSG}) we
obtain $P((xw^{-1})_{-\fS} \sqcup c \sqcup (xw^{-1})_{\fS}) =
P((xw^{\prime-1})_{-\fS} \sqcup c \sqcup (xw^{\prime-1})_{\fS})$, as
desired. \qed
\end{proof}

\begin{remark}\label{rem:BH}
The Propositions \ref{prop:induction}, \ref{prop:induction sn},
\ref{prop:domino lusztig} and \ref{prop:domino lusztig sn}
generalize \cite[Prop. 4.8]{BH} (which corresponds
to the asymptotic case).
\end{remark}


\subsection{Asymptotic case, quasi-split case}
We now prove Conjecture A for $r = 0, 1$ and 
$r \geq n-1$.

\begin{theorem}
\label{thm:BI} Conjecture~A is true for $r \geq n - 1$.
\end{theorem}

\begin{proof}
Let $r \geq n-1$, and let $D$ be a domino tableau with shape $\la
\in \cP_r(n)$.  The dominoes $\{\dom_i \mid i \in\{1,2,\ldots,n\}\}$
can be decomposed into the two disjoint collections $\cD_+ = \{\dom_i
\mid \mbox{$\dom_i$ is horizontal} \}$ and $\cD_- = \{\dom_i \mid
\mbox{$\dom_i$ is vertical} \}$ such that all the dominoes in $\cD_+$
lie strictly above and to the right of all the dominoes in $\cD_-$.
We call a tableau satisfying this property {\it segregrated}.  If
the collection of dominoes $\cD_+$ is left justified, and each domino
replaced by a single box, one obtains a usual Young tableau.
Similarly, if the dominoes $\cD_-$ are justified upwards and changed
into boxes, one obtains a usual Young tableau.

In other words, $D$ can be thought of as a union of two usual
tableau $D_+$ and $D_-$ so that the union of the values in $D_+$ and
in $D_-$ is the set $\{1,2,\ldots,n\}$.  To be consistent with the
remaining discussion we in fact define $D_-$ to be the {\it
conjugate} (reflection in the main diagonal) of the Young tableau
obtained from the dominoes $\cD_-$, as described above.

Domino insertion is compatible with this decomposition so that the
following diagram commutes:
\[\begin{picture}(300,60)
\put( 07,05){$[w_+,w_-]$}
\put( 60,06){\vector(1,0){75}}
\put( 70,10){\scriptsize{\text{RS insertion}}}
\put(150,05){$[P_+^r(w),Q_+^r(w)),(P^r_-(w),Q^r_-(w))]$}
\put( 20,50){$w$}
\put( 55,51){\vector(1,0){100}}
\put( 70,55){\scriptsize{\text{domino insertion}}}
\put(170,50){$(P^r(w),Q^r(w))$}
\put( 24,45){\vector(0,-1){28}}
\put(200,45){\vector(0,-1){28}}
\end{picture} \]
Here $w_+$ denotes the subword of $w$ consisting of positive letters
and $w_-$ denotes the subword consisting of negative letters, with
the minus signs removed.  In~\cite{BI}, it is shown for $r \geq n-1$
and $w, w' \in W_n$ that $w \sim_{\cL} w'$ if and only if $Q(w_+) =
Q(w'_+)$ and $Q(w_-) = Q(w'_-)$; similar results hold for
$\sim_{\cR}$ and $\sim_{\cLR}$. Since $Q(w_+) = Q^r_+(w)$ and
$Q(w_-) = Q^r_-(w)$, we have $w \sim_{\cL} w'$ if and only if
$Q^r(w) = Q^r(w')$, establishing Conjecture~A in this case
for left cells.  A similar argument works for right cells and
two-sided cells, using also the classification of two-sided cells
in~\cite{BI2}. \qed
\end{proof}

\begin{theorem}
Conjecture~A is true if $a=2b$ or if $3a=2b$.
\end{theorem}

\begin{proof}
In~\cite{Lusztig83}, Lusztig determined the left cells of $W_n$ with
parameters $b = (2r+1)a/2$ for $r \in \{0, 1\}$ as follows.  When $r
\in \{0,1\}$ we have $I = \emptyset$ in the notation of
Section~\ref{sec:lusztig}.  The equal parameter weight function $L$
on $\fS_{2n+r}$ restricts to the weight function $L_{b,a}$ on
$\iota_r(W_n)$, where $b = (2r+1)a/2$.
Lusztig~\cite[Theorem~11]{Lusztig83} shows that each left cell of
$W_n$ is the intersection of a left cell of $\fS_{2n+r}$ with
$\iota_r(W_n)$.  Thus $w \simeq_r w'$ in $W_n$ if and only if
$\iota_r(w) \simeq_\fS \iota_r(w')$ in $\fS_{2n+r}$.  When $r \in
\{0,1\}$ there is no need for the element $x \in \cb_0$ in
Theorem~\ref{thm:sn}, so one obtains Conjecture~A for
$r\in \{0,1\}$. \qed
\end{proof}


\subsection{Right descent sets}
If $r \geq 0$, let $S_n^{(r)} = \{s_1,s_2,\dots,s_{n-1}\} \cup
\{t_1,\dots,t_r\}$ (if $r \geq n$, then $S_n^{(r)}=S_n^{(n)}$).
If $w \in W_n$, let
$$\cR_n^{(r)}(w) = \{s \in S_n^{(r)}~|~\ell(ws) < \ell(w)\}$$
be the {\it extended right descent set} of $w$.
The following proposition is easy:

\begin{proposition}\label{right descent set}
Let $x$ and $y$ be two elements of $W_n$. Then:
\begin{itemize}
\item[(a)] If $x \simeq_r y$, then $\cR_n^{(r+1)}(x)=\cR_n^{(r+1)}(y)$.
\item[(b)] If $b > ra$ and if $x \sim_\cL y$, then $\cR_n^{(r+1)}(x)=
\cR_n^{(r+1)}(y)$.
\end{itemize}
\end{proposition}

\begin{proof}
If $r \geq n-1$, then statements (a) and (b) are equivalent by 
Theorem \ref{thm:BI}.
But, in this case, (b) has been proved in \cite[Prop. 4.5]{BH}.
So let us assume from now on that $r < n-1$. We shall prove (a) and
(b) together. Let us set
\begin{align*}
\cR_s(x)&=\{s \in \{s_1,\dots,s_{n-1}\}~|~\ell(ws) < \ell(w)\},\\
\cR_t^{(r)}(x)&=\{s \in \{t_1,\dots,t_r\}~|~\ell(ws) < \ell(w)\}.
\end{align*}
Then $\cR_n^{(r)}(x)=\cR_s(s) \cup \cR_t^{(r)}(x)$.

Write $x=ux'$ and $y=vy'$, with $u$, $v \in X_{r+1}^n$ and $x'$, $y' \in 
W_{r+1}$. Since $\ell(ux')=\ell(u)+\ell(x')$ (and similarly for $ux's$ for 
any $s \in W_r$), we have that $\cR_t^{(r)}(x)=\cR_t^{(r)}(x')$.
Similarly, $\cR_t^{(r)}(y)=\cR_t^{(r)}(y')$. But, if $x$ and $y$ satisfy
(a) or (b), then $\cR_t^{(r)}(x')=\cR_t^{(r)}(y')$:
indeed, this follows from the fact that (a) and (b) have been proved in the
asymptotic case and, in case (a), from Proposition \ref{prop:induction} and, 
in case (b), from \cite{geck induction}.

Now it remains to show that $\cR_s(x)=\cR_s(y)$ if $x$ and $y$ satisfy (a) 
or (b).  In case (b), this follows from \cite[Lemma 8.6]{Lusztig03}. So assume
now that $x \simeq_r y$. Write $x=u'\sigma$ and $y=v'\tau$, with $u$, 
$v \in X(n)$ and $\sigma$, $\tau \in \fS_n$. As in the previous case, we have
$\cR_s(x)=\cR_s(\sigma)$ and $\cR_s(y)=\cR_s(\tau)$. Moreover, by
Proposition \ref{prop:induction sn}, we have $\sigma \simeq_\fS \tau$.
It is well-known that it implies that $\cR_s(\sigma)=\cR_s(\tau)$. \qed
\end{proof}

\begin{remark}
Proposition~\ref{right descent set} (a) can also be deduced
from~\cite[Lemma~33]{ShWh} or~\cite[Lemma~9]{Lam1}.
\end{remark}


\subsection{Coplactic relations}
If $x$ and $y$ are two elements of $W_n$ such that $\ell(x) \leq \ell(y)$, 
then we write $x \smile_r y$ if there exists $s \in S_n^{(0)}$ and $s' \in 
S_n^{(r)}$ such that $y=sx$ and $\ell(s'x) < \ell(x) < \ell(y) < \ell(s'y)$. 
If $\ell(x) \geq \ell(y)$, then
we write $x \smile_r y$ if $y \smile_r x$. Let $\equiv_r$ denote
the reflexive and transitive closure of $\smile_r$.

\begin{remark}\label{coplactic easy}
(a) If $x \equiv_r y$, then $\ell_t(x)=\ell_t(y)$.

(b) If $r' \geq r$ and if $x \equiv_r y$, then $x \equiv_{r'} y$
(indeed, if $x \smile_{r} y$, then $x \smile_{r'} y$:
this just follows from the fact that $S_n^{(r)} \subset S_n^{(r')}$).
Moreover, the relations $\equiv_n$ and $\equiv_r$ are equal if $r \geq n$.

(c) If $r \geq n-1$, then $x \equiv_r y$ if and only if $x \equiv_{n-1} y$.
Let us prove this statement. By (b) above, we only need to show that,
if $x \smile_n y$, then $x \smile_{n-1} y$. For this, we may assume that
$\ell(y) > \ell(x)$. So there exists $i \in \{1,2,\dots,n-1\}$ and 
$s' \in S_n^{(n)}$ such that $y=s_i x$ and $\ell(s'x) < \ell(x) < \ell(s_ix) 
< \ell(s's_ix)$. If $s' \in S_n^{(n-1)}$ then we are done. So we may assume 
that $s'=t_n$. Therefore, the first inequality says that $x^{-1}(n) < 0$ and
the last inequality says that $x^{-1} s_i(n) > 0$. This implies that
$i=n-1$. Consequently, we have $x^{-1}(n) < 0$ and $x^{-1}(n-1) > 0$. But,
the middle inequality says that $x^{-1}(n-1) < x^{-1}(n)$, so we obtain
a contradiction with the fact that $s'=t_n$.

(d) If $b > ra$ and if $r \geq n-1$, then it follows from 
Theorem \ref{thm:BI}, from \cite[Prop. 3.8]{BI} and from (a), (b) and (c) 
above that the relations $\sim_\cL$, $\simeq_r$ and $\equiv_r$ coincide.
\end{remark}

\begin{proposition}\label{coplactic}
Let $x$ and $y$ be two elements of $W_n$ such that $x \equiv_r y$.
Then the following hold:
\begin{itemize}
\item[(i)] $x \simeq_r y$.
\item[(ii)] If $b > ra$, then $x \sim_\cL y$.
\end{itemize}
\end{proposition}

\begin{proof}
We may, and we will, assume that $x \smile_r y$. By symmetry, we may
also assume that $\ell(y) > \ell(x)$. We shall prove (i) and (ii)
together. There exists $s \in S_n^{(0)}$ and $s' \in S_n^{(r)}$ such
that $y=sx$ and $\ell(s'x) < \ell(x) < \ell(y) < \ell(s'y)$. Two
cases may occur:

\medskip

$\bullet$ If $s' \in S_n^{(0)}$, then write $x=x'u^{-1}$ and $y=y'v^{-1}$ with
$x'$, $y' \in \fS_n$ and $u$, $v \in X(n)$. Then $y'=sx'$, $u=v$ and
$\ell(s'x') < \ell(x') < \ell(y') < \ell(s'y')$. It is well-known that it 
implies that $Q(x')=Q(y')$ (Knuth relations), so $x' \simeq_\fS y'$
and $x' \sim_\cL y'$. Therefore, since moreover $u=v$,
it follows from Proposition \ref{prop:domino lusztig sn}
(resp. \cite[Prop. 9.13]{Lusztig03}) that $x \simeq_r y$
(resp. $x \sim_\cL y$).

\medskip

$\bullet$ If $s' \not\in S_n^{(0)}$ then we write $s=s_i$ and
$s'=t_j$. Then the relations $y=sx$ and $\ell(s'x) < \ell(x) <
\ell(y) < \ell(s'y)$ imply that $x^{-1}(j) < 0$ and
$x^{-1}s_i(j) > 0$. In particular, $s$ and $s'$ belong to
$W_{r+1}$. Now, write $x=x'u^{-1}$ and $y=y'v^{-1}$ with $x'$, $y'
\in W_{r+1}$ and $u$, $v \in X_{r+1}^n$. Then $y'=sx'$, $u=v$ and
$\ell(s'x') < \ell(x') < \ell(y') < \ell(s'y')$. By Remark
\ref{coplactic easy} (d), this implies that $x' \simeq_r y'$ and, if
$b > ra$, that $x' \sim_\cL y'$. Therefore, since moreover $u=v$, it
follows from Proposition \ref{prop:domino lusztig sn} (resp.
\cite[Prop. 9.13]{Lusztig03}) that $x \simeq_r y$ (resp. $x \sim_\cL
y$).\qed
\end{proof}

Even if we have both $\ell_t(x)=\ell_t(y)$ and $x \simeq_r y$ we do
not necessarily have $x \equiv_r y$.  For example, let $r = 0$, $n =
6$ and take $x = 5$ $6$ $1$ $4$ $2$ $-3$ and $y =$ $5$ $6$ $-1$ $4$
$3$ $2$.
%
%

\section{Cycles and Conjecture~B} \label{sec:cycles}

\subsection{Open and closed cycles}
We now describe a more refined combinatorial structure of domino tableaux 
introduced by Garfinkle~\cite{gar1}.  We will mostly follow the setup
of~\cite{pie1}.

Let $D$ be a domino tableau with shape $\la \in \cP_r(n)$.  We call a
square $(i,j) \in D$ {\it variable} if $i+j$ and $r$ have the same
parity, otherwise we call it {\it fixed}.  If the domino $\dom_i$
contains the square $(k,l)$ we write $D(k,l) = i$.

Now let $(k,l)$ be the fixed square of $\dom_i$.  Suppose that
$\dom_i$ occupies the squares $\{(k,l),(k+1,l)\}$ or
$\{(k,l-1),(k,l)\}$.  We define a new domino $\dom'_i$ by letting it
occupy the squares
\begin{enumerate}
\item
$\{(k,l),(k-1,l)\}$ if $i < D(k-1,l+1)$,
\item
$\{(k,l),(k,l+1)\}$ if $i > D(k-1,l+1)$.
\end{enumerate}
Otherwise $\dom_i$ occupies the squares $\{(k,l),(k,l+1)\}$ or
$\{(k-1,l),(k,l)\}$.  We define a new domino $\dom'_i$ by letting it
occupy the squares
\begin{enumerate}
\item
$\{(k,l),(k,l-1)\}$ if $i < D(k+1,l-1)$,
\item
$\{(k,l),(k+1,l)\}$ if $i > D(k+1,l-1)$.
\end{enumerate}

Now define the {\it cycle} $c = c(D,i)$ of $D$ through $i$ to be the
smallest union $c$ of dominoes satisfying that (i) $\dom_i \in c$
and (ii) $\dom_j \in c$ if $\dom_j \cap \dom'_k \neq \emptyset$
or $\dom'_j \cap \dom_k \neq \emptyset$ for some $\dom_k \in c$. If
$c$ is a cycle of $D$ we let $M(D,c)$ be the domino tableau obtained
from $D$ by replacing each domino $\dom_i \in c$ by $\dom'_i$.  We
call this procedure {\it moving through $c$}.

\begin{theorem}[\cite{gar1}]
\label{thm:cycle} Let $D$ be a domino tableau and $c$ a cycle of
$D$. Then $M(D,c)$ is a standard domino tableau.  Furthermore, if
$C$ is a set of cycles of $D$ then the tableau $M(D,C)$ obtained by
moving through each $c \in C$ is defined unambiguously.
\end{theorem}

We call a cycle $c$ {\it closed} if $M(D,c)$ has the same shape as
$D$; otherwise we call $c$ {\it open}.  Note that each (non-trivial)
cycle $c$ is in one of two positions, so that moving through is an
invertible operation.


\subsection{Evidence for Conjecture~B}
The notion of open and closed cycles allows us to state a
combinatorially more precise version of Conjecture~B.

\begin{conjA}
\label{conj4} Assume that $b=ra$ for some $r\geq 1$. Then the
following hold for any $w$, $w' \in W_n$:
\begin{itemize}
\item[(a)] $w \sim_\cL w'$ if and only if $Q^{r-1}(w) = M(Q^{r-1}(w'),C)$ 
for a set $C$ of {\itshape\bfseries open} cycles.
\item[(b)] $w \sim_\cR w'$ if and
only if $P^{r-1}(w)=M(P^{r-1}(w'),C)$ for a set $C$ of 
{\itshape\bfseries open} cycles.
\item[(c)] $w \sim_{\cLR} w'$ if and only if some tableau with shape
equal to $\sh(P^{r-1}(w)) = \sh(Q^{r-1}(w))$ can be obtained from a tableau 
with shape $\sh(P^{r-1}(w')) = \sh(Q^{r-1}(w'))$ by moving through a set of
{\itshape\bfseries open} cycles.
\end{itemize}
\end{conjA}

\begin{remark}
Each cycle $c$ of a domino tableau $D$ is in one of two positions,
and by Theorem~\ref{thm:cycle} they can be moved independently. Thus
Conjecture~D would imply that every left cell for the
parameters $b = ra$ would be a union of $2^d$ left cells for the
parameters $b = \frac{(2r-1)a}{2}$.  Here $d$ is equal to the number
of open cycles, which do not change the shape of the core, in one of
the $Q$-tableaux in Conjecture~D. This is consistent with the fact 
that, if $b=ra$ with $r \geq 1$, then the number of irreducible components 
of a constructible representation is a power of $2$ 
(see \cite[Chap. 22]{Lusztig03}; see also 
with~\cite[(12.1)]{Lusztig83} for the equal parameter case).
\end{remark}

We have the following theorem of Pietraho, obtained via a careful
study of the combinatorics of cycles.

\begin{theorem}[Pietraho \cite{pie2}]
Conjecture~B and Conjecture~D are equivalent.
\end{theorem}

Some special cases of Conjectures~B and~D are
known.  The case $b = a$ or $r = 1$ is known as the equal parameter
case and is closely connected with the classification of primitive
ideals of classical Lie algebras.

\begin{theorem}[Garfinkle \cite{gar2}]
Conjecture~D is true for $r = 1$.
\end{theorem}

The asymptotic case follows from~\cite{BI,BI2}.

\begin{theorem}
Conjecture~D holds for $r \geq n$.
\end{theorem}

\begin{proof}
Let $D$ be a domino tableau with shape $\la \in \cP_q(n)$ such that
$q \geq n -1$.  Then moving through any cycle of $D$ changes the
shape of the core of $D$.  Thus (for left cells) the condition
$Q^q(w) = M(Q^q(w'),C)$ in Conjecture~D is the same as the
condition $Q^q(w) = Q^q(w')$.  This agrees with the classification
given in~\cite{BI}.  A similar argument works for right cells and
two-sided cells, also using the classification in~\cite{BI2}. \qed
\end{proof}

\medskip
\noindent {\bf Acknowledgements.} 
T. L. was partially supported by NSF DMS--0600677.
Parts of the work presented here were done while all four authors 
enjoyed the hospitality of the Bernoulli center at the EPFL Lausanne 
(Switzerland), in the framework of the research program ``Group 
representation theory'' from January to June 2005.


\end{document}